\documentclass[12pt,twoside]{article} 
\usepackage{amsmath}
\usepackage{amssymb}

\newcommand{\Q}{\mbox{{\boldmath $Q$}}}

\title{%
\Large\bfseries
A note on the generalized $q$-Euler numbers(2)
}
\author{%
\vspace{5mm}
By \\
Kyoung-Ho {\sc Park}, Young-Hee {\sc Kim}, and Taekyun {\sc Kim}     }\vspace{15mm}
\def\Q{\mbox{\boldmath $Q$}}

\date{%
\begin{minipage}{14cm}%
\vspace{10mm}
\normalsize
{\bfseries {\small Abstract.}} {\small
Recently the new $q$-Euler numbers and polynomials related to Frobenius-Euler numbers and polynomials are constructed  by Kim (see[3]).
In this paper, we study the generalized $q$-Euler numbers and polynomials attached to $\chi$ related to the new $q$-Euler numbers
and polynomials which is constructed in [3].  Finally, we will derive some interesting congruence on the generalized $q$-Euler numbers and polynomials attached to $\chi$.
} 
\\
\end{minipage} \\
\begin{minipage}{14cm}
\normalsize
{\bfseries {\small 2000 Mathematics Subject Classification:}}
{\small
11S80, 11B68.\\
 }
 {\bfseries {\small Key Words and Phrases:}}
{\small
\vspace*{1cm}
 $q$-Euler numbers, $q$-Euler polynomials.
 }
 \\
\end{minipage}
}
\newcommand{\SECTION}[2]{%
  \vspace{5mm}\par\noindent $\S$#1. {\bf #2}~}

\newcommand{\THM}[1]{%
  \par\vspace{5mm}\par\noindent {\scshape #1.~}}

\newcommand{\REM}[1]{%
    \par\vspace{5mm}\par\noindent {\scshape #1.}}

\newcommand{\pn}{\par\noindent}

\newcommand{\bpr}{\bigskip\par}
\newcommand{\mpn}{\medskip\par\noindent}
\newcommand{\bpn}{\bigskip\par\noindent}
\begin{document}
\maketitle \pagestyle{myheadings} \markboth{{\footnotesize {\hfill
\rm A note on the generalized $q$-Euler numbers \hfill}}} {{ \small
{\hfill \rm Kyoung-Ho {\sc Park}, \rm Young-Hee {\sc Kim} and \rm Taekyun {\sc Kim}
  \hfill} }}
 \vspace* {-10mm}
 \vspace* {-25mm}
 \vspace* {-55mm}

\newcommand\QF[1]{\Q \hspace{-1mm} \left( \hspace{-1mm}\sqrt{#1} \right) }
\newcommand{\dsp}[1]{{\displaystyle #1}}
\newenvironment{THEOREM}[1]{%
    \vspace{5mm}\par\noindent{\bf Theorem. #1}~}{}
\newenvironment{REMARK}[1]{%
    \vspace{5mm}\par\noindent{\bf Remark #1}~}{}
\newenvironment{LEMMA}[1]{%
    \vspace{5mm}\par\noindent{\bf Lemma #1}~}{}
\newenvironment{PROPOSITION}[1]{%
    \vspace{5mm}\par\noindent{\bf Proposition #1}~}{}
\newenvironment{COROLLARY}[1]{%
    \vspace{5mm}\par\noindent{\bf Corollary. #1}~}{}
\newenvironment{DEFINITION}[1]{%
    \vspace{5mm}\par\noindent{\bf Definition #1}~}{}
\newenvironment{PROOF.}{%
    \vspace{5mm}\par\noindent{\bf Proof.}~}{%
    \hfill\hbox{\rule[-2pt]{3pt}{6pt}}%
    \par\vspace{5mm}}
\newenvironment{EXAMPLES}{%
    \vspace{5mm}\par\noindent{\sc Examples}.~}{}
\newenvironment{ACKNOWLEDGEMENTS}{%
    \vspace{5mm}\par\noindent{\bf Acknowledgements.}~}{}
\newcommand{\Mod}[1]{\,(\text{\mbox{\rm mod}}\;#1)}
\newcommand{\NUM}[1]{{[#1]}}
\baselineskip 6mm
\thispagestyle{empty}
\vspace{2cm}
\vspace{5cm}
\medskip
\par
\begin{center}\SECTION{1}{ Introduction}
\end{center}
Let $\mathbb{Z},\mathbb{R}$ and $\mathbb{C}$ denote the ring of integers, the field of real numbers and the complex number field.
and let $p$ be a fixed an odd prime number.Assume that $q$ is an indeterminate in $\Bbb C$ with $q \in \mathbb{C}$ with $|q| < 1.$
As the $q$-symbol $[x]_{q}$ , we denote $[x]_{q}=\frac{1-q^{x}}{1-q}$.
 Recently, $q$-Euler polynomials are defined as
 $$\frac{[2]_q}{qe^t+1}e^{xt}=\sum_{n=0}^{\infty} E_{n,q}(x)\frac{t^n}{n!}, \text{ for $|t+ \log q|<\pi,$ (see [3])}.$$
In the special case $x=0$, $E_{n,q}=E_{n,q}(0)$ are call the $n$-th $q$-Euler numbers (see [3]).
These $q$-Euler numbers and polynomials are closely relayed to Frobenius-Euler numbers and polynomials
and these numbers are studied by Simsek-Cangul-Ozden, Cenkci-Kurt and Can and several authors (see [1-2, 18-26]).
In this paper, we study the generalized $q$-Euler numbers and polynomials attached to $\chi$ related to the $q$-Euler numbers                 and polynomials, $E_{n,q}(x)$, which is constructed in [3].  Finally, we will derive some interesting congruence on
 the generalized $q$-Euler numbers and
polynomials attached to $\chi$.                                                                                                                   
\begin{center}
\SECTION{2}{ Congruence for $q$-Euler numbers and polynomials }
\end{center}
The ordinary Euler polynomials are defined as
\begin{equation*}
e^{xt}\frac{2}{e^{t}+1}=e^{E(x)t}=\sum_{n=0}^{\infty}E_{n}(x)\frac{t^{n}}{n!}, \quad ({\rm see \ [1-5]}),
\end{equation*}
where we use the technical method notation by replacing $E^{n}(x)$ by $E_{n}(x)(n \geq 0)$, symbolically (see [1-2]).
Let us consider the generating function of $q$-Euler polynomials $E_{n,q}(x)$ as follows:
\begin{equation*}
\tag{1}
F_{q}(x,t)=\frac{[2]_{q}}{qe^t +1}e^{xt}=\sum_{n=0}^{\infty}E_{n,q}(x)\frac{t^{n}}{n!},
\end{equation*}
and we also note that
$$\sum_{n=0}^{\infty}E_{n,q}(x)\frac{t^n}{n!}=\frac{[2]_q}{qe^t+1}e^{xt}=\frac{1-(-q^{-1})}{e^t-(-q^{-1})}=\sum_{n=0}^{\infty}H_n(-q^{-1},x)\frac{t^n}{n!},$$
where $H_n(-q^{-1},x) $ are called the $n$-th Frobenius-Euler polynomials (see [3]).
From (1), we note that
\begin{equation*}
\tag{2}
\lim_{q\rightarrow 1}F_{q}(x,t)=\frac{2}{e^{t}+1}e^{xt}=\sum_{n=0}^{\infty}E_n(x)\frac{t^n}{n!}.
\end{equation*}
By (1) and (2), we see that
\begin{equation*}
\lim_{q\rightarrow 1}E_{n,q}(x)=E_{n}(x).
\end{equation*}
In (1), it is easy to show that

$$\sum_{n=0}^{\infty}E_{n,q}(x)\frac{t^{n}}{n!}=F_{q}(x,t)=\frac{[2]_q}{qe^{t}+1}e^{xt}
=\sum_{n=0}^{\infty} \Big( \sum_{l=0}^n \binom{n}{l}E_{l, q} x^{n-l} \Big)\frac{t^{n}}{n!}.$$
By comparing the coefficients on the both sides, we have
\begin{equation*}
\tag{3}
E_{n,q}(x)=\sum_{l=0}^n \binom{n}{l}E_{l, q} x^{n-l}, \quad \text{where $E_{l,q}$ are the $l$-th $q$-Euler numbers}.
\end{equation*}
Let $\chi$ be the Dirichlet's character with conductor $d \equiv 1\Mod{2}$. Then we define generating function of the generalized
$q$-Euler numbers attached to $\chi$, $E_{n,\chi,q}$ as follows:
\begin{equation*}
\tag{4}
F_{q,\chi}(t)=\frac{[2]_{q}\sum_{l=0}^{d-1}\chi(l)q^l(-1)^{l}e^{lt}}{q^de^{dt}+1}=\sum_{n=0}^{\infty}E_{n,\chi,q}\frac{t^{n}}{n!}.
\end{equation*}
From (4), we note that
\begin{equation*}
\tag{5}
\lim_{q\rightarrow 1}F_{q,\chi}(t)=\frac{2\sum_{a=0}^{d-1}\chi(a)(-1)^{a}e^{at}}{e^{dt}+1}=\sum_{n=0}^{\infty}E_{n,\chi}\frac{t^{n}}{n!},
\end{equation*}
where $E_{n,\chi}$ are the $n$-th ordinary Euler numbers attached to $\chi$. By (4) and (5), we see that
\begin{equation*}
\lim_{q\rightarrow 1}E_{n,\chi,q}=E_{n,\chi}.
\end{equation*}
From (5), we can also derive
\begin{equation*}
\tag{6}
\begin{split}
\sum_{n=0}^{\infty}E_{n,\chi, q}\frac{t^{n}}{n!}&=F_{q,\chi}(t)=[2]_{q}\sum_{k=0}^{\infty}\chi(k)(-q)^{k}e^{kt}\\
&=\sum_{n=0}^{\infty}\Big([2]_{q}\sum_{k=0}^{\infty}\chi(k)(-q)^{k}k^{n} \Big)\frac{t^{n}}{n!}\\
&=\sum_{n=0}^{\infty}\Big(d^n\sum_{a=0}^{d-1}(-q)^{a}\chi(a)
E_{n,q^d}(\frac{a}{d}) \Big)\frac{t^{n}}{n!}.
\end{split}
\end{equation*}
By comparing the coefficients on the both sides of (6), we have
\begin{equation*}
\tag{7}
\begin{split}
E_{n,\chi, q}&=[2]_{q}\sum_{k=0}^{\infty}\chi(k)(-q)^{k}k^{n}=d^n\sum_{a=0}^{d-1}(-q)^{a}\chi(a)E_{n,q^d}(\frac{a}{d}) .
\end{split}
\end{equation*}
Finally, we define the generating function of the generalized
$q$-Euler polynomials attached to $\chi$, $E_{n,\chi,q}(x)$ as
follows:
\begin{equation*}
\tag{8}
F_{q,\chi}(x,t)=\sum_{n=0}^{\infty}E_{n,\chi,q}(x)\frac{t^{n}}{n!}=[2]_{q}\sum_{k=0}^{\infty}\chi(k)(-q)^{k}e^{(x+k)t}.
\end{equation*}
By (8), we easily see that
\begin{equation*}
\tag{9}
\begin{split}
\sum_{n=0}^{\infty}E_{n,\chi,q}(x)\frac{t^{n}}{n!}
&=F_{q,\chi}(x,t)=[2]_{q}\sum_{k=0}^{\infty}\chi(k)(-q)^{k}e^{(x+k)t}\\
&=\sum_{n=0}^{\infty}\Big( [2]_{q}\sum_{k=0}^{\infty}\chi(k)(-q)^{k}(x+k)^{n}\Big)\frac{t^{n}}{n!}\\
&=\sum_{n=0}^{\infty}\Big(d^n\sum_{a=0}^{d-1}(-q)^{a}\chi(a) E_{n,q^d}(\frac{a+x}{d}) \Big)\frac{t^{n}}{n!}.
\end{split}
\end{equation*}
Thus, we have
\begin{equation*}
\tag{10}
\begin{split}
E_{n,\chi,q}(x)
&=d^n\sum_{a=0}^{d-1}(-q)^{a}\chi(a) E_{n,q^d}(\frac{a+x}{d}) =\sum_{\ell=0}^{n}\binom{n}{\ell}x^{n-\ell}E_{\ell, \chi, q}
=[2]_{q}\sum_{k=0}^{\infty}\chi(k)(-q)^{k}(x+k)^{n}.
\end{split}
\end{equation*}
 Let $d\in \mathbb{N}$ with $d \equiv 1\Mod{2}$. Then, we see that
\begin{equation*}
\tag{11}
\begin{split}
q^dF_{q,\chi}(d,t)+F_{q,\chi}(t)
&=[2]_{q}\sum_{k=0}^{\infty}\chi(k)(-q)^{k}e^{ (d+k)t}
+ [2]_{q}\sum_{k=0}^{\infty}\chi(k)(-q)^{k}e^{ kt}\\
&=[2]_{q}\sum_{k=0}^{d-1}\chi(k)(-q)^{k}e^{ kt}.
\end{split}
\end{equation*}
From (11), we have
\begin{equation*}
\sum_{n=0}^{\infty}\Big(q^dE_{n,\chi,q}(d)+E_{n,\chi,q}\Big)\frac{t^{n}}{n!}
=\sum_{n=0}^{\infty}\Big\{[2]_{q}\sum_{k=0}^{d-1}\chi(k)(-q)^{k}k^{n}\Big\} \frac{t^{n}}{n!}.
\end{equation*}
Therefore, we obtain the following theorem.
\medskip
{\THM {Theorem {1}}} {\it For $q \in \mathbb{C}$ with $|q|<1, n\in \mathbb{Z}_{+}$ and $d \in \mathbb{N}$ with $d\equiv 1\Mod{2}$, we have
\begin{equation*}
q^dE_{n,\chi,q}(d)+E_{n,\chi,q}=[2]_{q}\sum_{k=0}^{d-1}\chi(k)(-q)^{k}k^{n}.
\end{equation*}
\rm}
Let $p$ be a positive odd integer and let $N \in \mathbb{N}$. Then we have
\begin{equation*}
\begin{split}
[2]_{q}\sum_{a=0}^{dp^{N}-1}\chi(a)(-q)^{a}a^{n}
&=q^{dp^N}E_{n,\chi, q}(dp^{N})+E_{n,\chi,q}\\
&=q^{dp^N}\sum_{j=0}^{n}\binom{n}{j}(dp^{N})^{j}E_{n-j,\chi,q}+E_{n,\chi,q}\\
&=q^{dp^N}\sum_{j=1}^{n}\binom{n}{j}(dp^{N})^{j}E_{n-j,\chi,q}+(q^{dp^{N}}+1)E_{n,\chi,q}\\
&\equiv 2E_{n,\chi,q} \Mod{dp^{N}},
\end{split}
\end{equation*}
because $q^{ndp^{N}}\equiv 1\Mod{dp^{N}}$.
Therefore, we obtain the following theorem.
{\THM {Theorem {2}}} {\it Let $p$ be a positive odd integer and $q \in \mathbb{C}$ with $|q|<1$ and $(q-1, dp)=1$. For
 $d \in \mathbb{N}$ with $d\equiv 1\Mod{2}$, we have
\begin{equation*}
[2]_{q}\sum_{a=0}^{dp^{N}-1}\chi(a)(-q)^{a}a^{n}\equiv 2E_{n,\chi,q} \Mod{dp^{N}}.
\end{equation*}
\rm}
\medskip
{\REM {Remark}} {\it Define
\begin{equation*}
L_{E,q}(s,\chi|x)=[2]_{q}\sum_{n=0}^{\infty}\frac{(-q)^{n}\chi(n)}{(n+x)^{s}},
\end{equation*}
where $s\in \mathbb{C}$, and $x \neq 0, -1, -2, \cdots .$ For $k \in \mathbb{Z}_{+}$, we have $L_{E,q}(-k, \chi |x)=E_{k,\chi,q}(x)$.
\rm}
\medskip
\bpr

  ACKNOWLEDGEMENTS. The present research has been conducted by the
research grant of the Kwangwoon University in 2009.

\bigskip\bigskip
\begin{center}\begin{large}
{\sc References}
\end{large}\end{center}
\par
\begin{enumerate}
\item[{[1]}] R. P. Agarwal, C. S. Ryoo, {\it Numerical Computations of the roots of the generalized twisted $q$-Bernoulli polynomials},
Neural parallel Sci. Comput, {\bf 15} (2007), 193-206.

\item[{[2]}] I. N. Cangul, H. Odzen and Y. Simsek, {\it Generating Functions of the $(h,q)$-extension of Euler polynomials and numbers}, Acta Math. Hungar. {\bf 120} (2008), 281-299.

\item[{[3]}] T. Kim, {\it An invariant $p$-adic $q$-integral on $\Bbb Z_p$},
Appl. Math. Letters {\bf 21} (2008), 105-108.

\item[{[4]}] T. Kim, {\it $q$-extension of the Euler formula and trigonometric functions}, Russ. J. Math. Phys. {\bf 14 (3)} (2007),275-278.

\item[{[5]}] T. Kim, {\it Symmetry $p$-adic invariant integral on $\mathbb{Z}_{p}$ for Bernoulli and Euler polynomials},
J. Difference Equ. Appl. {\bf 14} (2008), 1267-1277.

\item[{[6]}] H. M. Srivastava, T. Kim, and Y. Simsek, {\it $q$-Bernoulli numbers and polynomials associated with multiple $q$-zeta functions and basic $L$-series}, Russ. J. Math. Phys. {\bf 12(2)} (2005), 241-268.

\item[{[7]}] T. Kim, {\it Note on $q$-Genocchi numbers and polynomials}, Adv. Stud. Contemp. Math. {\bf 17} (2008), 9-15.

\item[{[8]}] T. Kim, {\it The modified $q$-Euler numbers and polynomials}, Adv. Stud. Contemp.
Math. {\bf 16} (2008), 161-170.

\item[{[9]}] T. Kim, {\it On a $q$-analogue of the $p$-adic log gamma functions and related integrals}, J. Number Theory,
{\bf 76} (1999), 320-329.

\item[{[10]}] T. Kim, {\it $q$-Volkenborn Integration}, Russ. J. Math. Phys. {\bf 9} (2002), 288-299.

\item[{[11]}] T. Kim, {\it $q$-Bernoulli numbers and polynomials associated with Gaussian binomial coefficients},
 Russ. J. Math. Phys. {\bf 15} (2008), 51-57.

\item[{[12]}] T. Kim, J. Y. Choi, and J. Y. Sug, {\it Extended $q$-Euler numbers and polynomials associated with fermionic $p$-adic $q$-integral on $\mathbb{Z}_{p}$}, Russ. J. Math. Phys. {\bf 14} (2007), 160-163.

\item[{[13]}] T. Kim, {\it On the von Staudt-Clausen's Theorem for the $q$-Euler numbers}, Russ. J. Math. Phys.
{\bf 16 (3)} (2009).

\item[{[14]}] T. Kim, {\it $q$-generalized Euler numbers and polynomials}, Russ. J. Math. Phys.{\bf 13(3)}(2006), 293-298.

\item[{[15]}] T. Kim, {\it Multiple $p$-adic $L$-function} Russ. J. Math. Phys. {\bf 13(2)} (2006), 151-157.

\item[{[16]}] T. Kim, {\it Power series and asymptotic series associated with the $q$-analog of the two-variable $p$-adic $L$-function}, Russ. J. Math. Phys. {\bf 12(2)} (2005), 186-196.

\item[{[17]}] T. Kim, {\it Analytic continuation of multiple $q$-zeta functions and their values at negative integers}, Russ. J. Math. Phys. {\bf 11(1)} (2004), 71-76.

\item[{[18]}] T. Kim {\it On Euler-Barnes multiple zeta functions}, Russ.J. Math. Phys. {\bf 10(3)} (2003), 261-267.

\item[{[19]}] B. A. Kupershmidt, {\it Reflection symmetries of $q$-Bernoulli polynomials }, J. Nonlinear Math. Phys. {\bf 12} (2005), 412-422.

\item[{[20]}] H. Ozden, Y. Simsek, {\it Interpolation functions of the $(h,q)$-extension of twisted Euler numbers}, Comput. Math. Appl. {\bf 56} (2008), 898-908.

\item[{[21]}] H. Ozden, Y. Simsek, {\it A new extension of $q$-Euler numbers and polynomials related to their interpolation functions}, Appl. Math. Lett. {\bf 21} (2008), 934-939.

\item[{[22]}] M. Schork, {\it Ward's calculus sequence; $q$-calculus and their limit $q \rightarrow -1$}, Adv. Stud. Contemp. Math. {\bf 13} (2006), 131-141.

\item[{[23]}] Y. Simsek, {\it Generating functions of the twisted Bernoulli numbers and polynomials associated with their interpolation functions}, Adv. Stud. Contemp. Math. {\bf 16} (2008), 251-278.

\item[{[24]}] Y. Simsek, {\it $p$-adic twisted $q\text{-}L$-functions related to generalized twisted Bernoulli numbers}, Russ. J. Math. Phys. {\bf 13(3)} (2006), 340-348.

\item[{[25]}] Y. Simsek, O. Yurekli, and V. Kurt, {\it On interpolation functions of the twisted generalized Frobenius-Euler numbers}, Adv. Stud. Contemp. Math. {\bf 15(2)} (2007),187-194.

\item[{[26]}] S. H. Rim, K. H. Park, and E.J. Moon, {\it On Genocchi numbers and polynomials}, Abstract and Applied Analysis, Article ID {\bf 898471}, (2008) 7 pages, dio:10.1155/2008/898471.
 \item[{[27]}]  T. Kim, Y.H. Kim, K. W. Hwang {\it On the $q$-extensions of the Bernoulli and Euler numbers,
     related identities and lerch zeta function}, Proc. Jangjeon Math. Soc. {\bf 12}
          (2009), no. 1, 1798--1804.

\end{enumerate}

\newpage
\par\noindent
\mpn { \bpn {\small Kyoung-Ho {\sc Pak}  \mpn Department of Mathematics, \pn Sogang University, Seoul 121-741, S.
Korea \pn {\it E-mail:}\ {\sf sagamath@yahoo.co.kr}
}
 \mpn { \bpn {\small Young-Hee {\sc Kim}  \mpn Division of General
Education-Mathematics, \pn Kwangwoon University, Seoul 139-701, S.
Korea \pn {\it E-mail:}\ {\sf yhkim@kw.ac.kr}
}
  \mpn { \bpn {\small Taekyun {\sc Kim}  \mpn Division of General
Education-Mathematics, \pn Kwangwoon University, Seoul 139-701, S.
Korea \pn {\it E-mail:}\ {\sf tkkim@kw.ac.kr}
}

\end{document}